\documentclass[12pt]{amsart}

\makeatletter
\@namedef{ver@rerunfilecheck.sty}{2026/07/05}
\makeatother
\newcommand{\RerunFileCheck}[4]{}

\usepackage{amsfonts, amsthm, amsmath, dsfont}
\usepackage{amssymb}
\usepackage{amscd}
\usepackage{lmodern}
\usepackage[latin2]{inputenc}
\usepackage{t1enc}
\usepackage[mathscr]{eucal}
\usepackage{mathrsfs}
\usepackage{mathtools} 

\usepackage{mathabx} 

\usepackage{graphicx}
\usepackage{graphics}
\usepackage{rotating}
\usepackage{pict2e}
\usepackage{epic}
\usepackage{framed}
\usepackage{booktabs}
\usepackage{makecell}
\usepackage{comment}
\usepackage{indentfirst}
\usepackage{enumerate}

\usepackage{tikz}
\usepackage{mathdots}
\usepackage{pgfplots}
\usetikzlibrary{positioning, arrows.meta, tikzmark}

\newcolumntype{V}{!{\vrule width 2pt}}

\usepackage[pagebackref]{hyperref}
\hypersetup{colorlinks=true}
\usepackage{bookmark} 

\usepackage{cite}
\usepackage{color}
\usepackage{refcount}

\allowdisplaybreaks
\numberwithin{equation}{section}

\theoremstyle{plain}
\newtheorem{theorem}{Theorem}[section]

\newtheorem{lemma}[theorem]{Lemma}
\newtheorem{definition}[theorem]{Definition}

\newtheorem{fact}[theorem]{Fact}

\def\eoc{\mathrm{eoc}}
\def\pom{\mathrm{pom}}
\def\emac{\mathrm{emac}}
\def\smic{\mathrm{smic}}
\def\smac{\mathrm{smac}}
\def\N{\mathbb{N}}
\def\fS{\mathfrak{S}}
\def\fT{\mathfrak{T}}
\def\AA{\mathfrak{AA}}
\def\rc{\mathrm{rc}}
\def\trc{\widetilde{\mathrm{rc}}}

\newcommand{\comp}{\mathbin{\scalebox{0.75}{$\circ$}}}

\def\multiset#1#2{\ensuremath{\left(\kern-.3em\left(\genfrac{}{}{0pt}{}{#1}{#2}\right)\kern-.3em\right)}}
\newcommand{\bintree}[1]{%
  \begin{tikzpicture}[
    baseline=(current bounding box.center),
    x=0.42cm, y=0.50cm,
    every node/.style={inner sep=1.5pt, font=\footnotesize}
  ]
    #1
  \end{tikzpicture}%
}

\begin{document}
	
	\title[Symmetric joint distribution on strict binary trees]{A combinatorial proof of a symmetric joint distribution on strict binary trees}
	
	\author[S. Fu]{Shishuo Fu}
	\address[S. Fu]{College of Mathematics and Statistics \& Center for Discrete Mathematics, Chongqing University, Chongqing 401331, China}
	\email{fsshuo@cqu.edu.cn}
	\dedicatory{Dedicated to Dominique Foata, master of bijective combinatorics.}

	\date{\today}

	\begin{abstract} 
	Nonplanar strict binary trees on $2n+1$ vertices are enumerated by the reduced tangent numbers $T_{2n+1}/2^n$, whose entries are refined by Poupard's triangle. Poupard defined the end of the minimal chain ($\mathrm{eoc}$) and the parent of the maximum leaf ($\mathrm{pom}$) on these trees, proving that $\mathrm{eoc}$ and $\mathrm{pom}+1$ share this triangular distribution. Foata and Han subsequently showed via tree calculus and partial difference equations that the joint distribution of $(\mathrm{eoc}, 2n+1-\mathrm{pom})$ is symmetric. By introducing an auxiliary tree statistic dual to $\mathrm{eoc}$, we give a bijective proof of this symmetry. Our proof composes a recursive bijection on strict binary trees, a bijection to alternating Andr\'e permutations, and the reverse-complement involution on permutations.
	\end{abstract}

	
	\maketitle

	\section{Introduction}
	The \emph{tangent numbers} (see, e.g.~\cite[p.~258--259]{Comtet1974}, \cite[A000182]{OEIS}) refer to the coefficients $T_{2n+1}$ in the Taylor expansion
	\begin{align}\label{gf:tangent}
	\tan x &= \sum_{n\ge 0}T_{2n+1}\frac{x^{2n+1}}{(2n+1)!}=1\frac{x}{1!}+2\frac{x^3}{3!}+16\frac{x^5}{5!}+272\frac{x^7}{7!}+7936\frac{x^9}{9!}+\cdots.
	\end{align}
	They bear rich combinatorial and arithmetical meanings and properties. For instance, a classical result of Andr\'e~\cite{Andre1879} in 1879 states that $T_{2n+1}$ counts \emph{alternating permutations} of length $2n+1$; see Stanley's survey \cite{Stanley2010} and the references therein. It is also well known that the set $\fT_{2n+1}$ of unordered strict binary trees (see Definition~\ref{defn:strict binary tree}) is enumerated by the \emph{reduced tangent number} $T_{2n+1}/2^n$; see \cite{FoataHan2014} and \cite[A002105]{OEIS}.

	In 1989, Christiane Poupard~\cite{Poupard1989} initiated the study on the system of finite difference equations
	\begin{align}\label{eq:diff eq}
	\Delta^2 f_n(m)+2f_{n-1}(m)=0,\quad 1\le m\le 2n-1,
	\end{align}
	where $\Delta$ stands for the finite difference operator
	\begin{align*}
	\Delta f_n(m) := f_n(m+1)-f_n(m).
	\end{align*}
	Together with the initial values $f_0(1)=1$, $f_n(1)=0$, and $f_n(2)=\sum_m f_{n-1}(m)$ for $n\ge 1$, the system~\eqref{eq:diff eq} produces a triangular array that was later named after Poupard; see Table~\ref{tab:poupard-triangle} for its first five rows. In the same paper~\cite{Poupard1989}, Poupard introduced two statistics ($\eoc$ and $\pom$, whose definitions are deferred to Section~\ref{sec:preliminaries}) defined over $\fT_{2n+1}$ to combinatorially account for each entry $f_n(m)$, thereby showing that the so-called \emph{Poupard's triangle} $f:=(f_n(m))$ refines the reduced tangent number in the following sense
	\[
		\sum_{m=1}^{2n+1} f_n(m) = \frac{T_{2n+1}}{2^n}.
	\]
	\begin{theorem}[Poupard~\cite{Poupard1989}]\label{thm:Poupard}
	For all $1\le m\le 2n+1$, we have
	\begin{align}\label{single Poupard}
	\#\{T\in\fT_{2n+1}:\eoc(T)=m\} = \#\{T\in\fT_{2n+1}:\pom(T)=m-1\}=f_n(m).
	\end{align}
	\end{theorem}
	Poupard's original proof of Theorem~\ref{thm:Poupard} makes use of a clever finite difference analysis. We note that an explicit bijection inspired from the classical ``jeu de taquin'' was given by Han~\cite{Han2013} to show the first equality of \eqref{single Poupard}.

	\begin{table}[htbp]
	\centering
	\caption{Poupard's triangle.}
	\label{tab:poupard-triangle}
	\renewcommand{\arraystretch}{1.25}
	$\begin{array}{r*{9}{c}}
	\hline\noalign{\medskip}
	 & & & & & f_0(1) & & & & \\
	 & & & & f_1(1) & f_1(2) & f_1(3) & & & \\
	f = & & & f_2(1) & f_2(2) & f_2(3) & f_2(4) & f_2(5) & & \\
	 & & f_3(1) & f_3(2) & f_3(3) & f_3(4) & f_3(5) & f_3(6) & f_3(7) & \\
	 & f_4(1) & f_4(2) & f_4(3) & f_4(4) & f_4(5) & f_4(6) & f_4(7) & f_4(8) & f_4(9) \\[1.5ex]
	 & & & & & 1 & & & & \\
	 & & & & 0 & 1 & 0 & & & \\
	= & & & 0 & 1 & 2 & 1 & 0 & & \\
	 & & 0 & 4 & 8 & 10 & 8 & 4 & 0 & \\
	 & 0 & 34 & 68 & 94 & 104 & 94 & 68 & 34 & 0 \\
	\noalign{\medskip}\hline
	\end{array}$
	\end{table}

	Now let us consider the joint distribution of both $\eoc$ and $\pom$ to define
	\[
		f_n(m,k) := \#\{t\in\fT_{2n+1}:\eoc(t)=m,~\pom(t)=k\},
	\]
	and introduce the sequence of matrices $\{M_n\}_{n\ge 1}$, where $M_n:=(f_n(m,k))_{1\le m,k\le 2n}$. The first three matrices are shown below, from which the anti-diagonal symmetry between the entries is evident. This stronger joint symmetric distribution was shown by Foata and Han~\cite{FoataHan2014} to hold in general.
	\[
	M_1 := 
	\begin{pmatrix} 
	0 & 0 \\ 
	1 & 0 
	\end{pmatrix},
	\quad
	M_2 = 
	\begin{pmatrix} 
	0 & 0 & 0 & 0 \\ 
	0 & 0 & 1 & 0 \\ 
	1 & 1 & 0 & 0 \\ 
	0 & 1 & 0 & 0 
	\end{pmatrix},
	\quad
	M_3 = 
	\begin{pmatrix} 
	0 & 0 & 0 & 0 & 0 & 0 \\ 
	0 & 0 & 1 & 2 & 1 & 0 \\ 
	1 & 1 & 0 & 4 & 2 & 0 \\ 
	2 & 3 & 4 & 0 & 1 & 0 \\ 
	1 & 3 & 3 & 1 & 0 & 0 \\ 
	0 & 1 & 2 & 1 & 0 & 0 
	\end{pmatrix}.
	\]
	\begin{theorem}[Foata-Han~{\cite[Theorem 1.2]{FoataHan2014}}]\label{thm:joint symm}
	For all $n\ge 1$, the matrices $M_n$ are symmetric with respect to their anti-diagonals:
	\begin{align}\label{eq:f-symm}
	f_n(m,k) &= f_n(2n+1-k,2n+1-m).
	\end{align}
	\end{theorem}

	Foata and Han's proof of Theorem~\ref{thm:joint symm} is indirect and relies on an in-depth analysis of the joint distribution matrices $M_n$. Specifically, they establish two dual systems of partial finite difference equations (the Delta and Gamma sequences) on the one hand, and develop an intricate ``tree calculus'' involving extensive case-by-case subtree decompositions on the other hand. As Han himself later remarked\footnote{Personal communication in Feb 2026.}, verifying these multi-step cancellations is quite involved and delicate. The purpose of this note is to provide a short, direct combinatorial proof of this symmetry. 

	A key ingredient of our approach is a new tree statistic called ``$\emac$'' (see Definition~\ref{def:three stats}) that is in a sense dual to $\eoc$ and it provides an alternative interpretation of $f_n(m,k)$.

	\begin{theorem}\label{thm:eoc to emac}
	There exists a bijection $\Phi: \fT_{2n+1}\to \fT_{2n+1}$ satisfying that for any $t\in\fT_{2n+1}$,
	\begin{align}\label{pom=pom}
	\pom(\Phi(t)) &= \pom(t),\\
	\emac(\Phi(t)) &= \eoc(t). \label{emac=eoc}
	\end{align}
	Consequently, we have 
	\begin{align}\label{eq:f as emac and pom}
	f_n(m,k) = \#\{t\in\fT_{2n+1}: \emac(t)=m,~\pom(t)=k\}.
	\end{align}
	\end{theorem}

	The new interpretation in \eqref{eq:f as emac and pom} is further transformed to a bivariate refinement of alternating Andr\'e permutations, another model that is counted by the reduced tangent number. Let $\AA_{n}$ denote the set of $n$-permutations that are simultaneously alternating (down-up) and Andr\'e. It turns out that this final interpretation over $\AA_n$ makes the anti-diagonal symmetry \eqref{eq:f-symm} transparent.

	\begin{theorem}\label{thm:tree to perm}
	There exists a bijection $\Psi:\fT_{2n+1}\to \AA_{2n+1}$ such that for every $t\in\fT_{2n+1}$ and $\pi:=\Psi(t)$, we have
	\begin{align}\label{stats:tree to perm}
	\emac(t) = \pi_1,\quad \pom(t) = \pi_{2n}.
	\end{align}
	Moreover, for any $1\le m,k\le 2n$,
	\begin{align}\label{eq:reverse-complement}
	\#\{\pi\in\AA_{2n+1}: \pi_1=m,~\pi_{2n}=k\} = \#\{\pi\in\AA_{2n+1}: \pi_1=2n+1-k,~\pi_{2n}=2n+1-m\}.
	\end{align}
	\end{theorem}

	Combining Theorems~\ref{thm:eoc to emac} and~\ref{thm:tree to perm} yields a direct bijective proof of Theorem~\ref{thm:joint symm}. We elaborate on this in the next section, where we introduce the necessary definitions and prove Theorem~\ref{thm:tree to perm}. The proof of Theorem~\ref{thm:eoc to emac} is given in Section~\ref{sec:eoc-emac}, which ends with a complete worked example illustrating all of the bijections.

	\section{Preliminaries}\label{sec:preliminaries}

	A \emph{binary tree} on the label set $L\subset \N$ is defined recursively by having a root vertex $r$ and a left and right subtree of $r$ which are themselves binary trees or are empty. All trees considered in this note are binary and \emph{increasingly labeled} by the set $L$, i.e., along the path from the root $r=\min(L)$ to any vertex $v\in L$ the labels are increasing. Equivalently, the root of any subtree receives the smallest label among all vertices on this subtree. When no confusion is caused, we speak of a vertex and its label interchangeably. We will mainly focus on the following kind of binary trees.

	\begin{definition}\label{defn:strict binary tree}
	A binary tree is called \emph{strict} (or complete) if every vertex has out-degree $0$ (a leaf) or $2$ (an internal node with two children). A strict binary tree must have an odd number, say $2n+1$, of vertices, consisting of $n$ internal nodes and $n+1$ leaves. A binary tree is called \emph{nonplanar} (or unordered) if the order of the two children of any internal node is irrelevant. We denote by $\fT_{2n+1}$ the set of nonplanar strict binary trees on $[2n+1]:=\{1,2,\ldots,2n+1\}$. 
	\end{definition}

	For a labeled tree $t$ with root labeled by $r$, we use $\min(t)$ and $\max(t)$ to denote the minimal and maximal label on $t$, respectively. Since all trees considered here are increasingly labeled, we always have $\min(t)=r$. Let $a$ be the label of an internal node of $t$ having two children labeled $b$ and $c$, and suppose $t_b$ and $t_c$ are the subtrees rooted at $b$ and $c$, respectively. We define 
	\begin{align}
		\smic(a) &:= \min\{\min(t_b),\min(t_c)\}=\min\{b,c\},\label{def:min}\\ 
		\smac(a) &:= \begin{cases}
		b & \text{if $\max(t_b)<\max(t_c)$,}\\
		c & \text{if $\max(t_b)>\max(t_c)$.}
		\end{cases}
		\label{def:max}
	\end{align}

	Now we introduce three statistics on the set of nonplanar strict binary trees, the second of which appears to be new and plays a key role in our bijective approach to Theorem~\ref{thm:joint symm}. 

	\begin{definition}\label{def:three stats}
	For $n\ge 0$, let $t$ be a nonplanar strict binary tree with label set $L$, and set $m:=\min(t)$ and $M:=\max(t)$.
	\begin{enumerate}
		\item The \emph{(smaller) minimal chain} of $t$ is defined to be the sequence $a_1\to a_2\to\cdots\to a_j$ that satisfies the following conditions:
		\begin{enumerate}[(i)]
			\item $a_1=m$ is the label of the root of $t$ and $a_j$ is the label of a leaf;
			\item for $i=1,2,\ldots,j-1$, the $i$-th term $a_i$ is the label of an internal node and $a_{i+1}=\smic(a_i)$.
		\end{enumerate}
		Define the \emph{end of the minimal chain} of $t$ to be: $\eoc(t):=a_j$.
		\item The \emph{(smaller) maximal chain} of $t$ is defined to be the sequence $b_1\to b_2\to\cdots\to b_k$ that satisfies the following conditions:
		\begin{enumerate}[(i)]
			\item $b_1=m$ is the label of the root of $t$ and $b_k$ is the label of a leaf;
			\item for $i=1,2,\ldots,k-1$, the $i$-th term $b_i$ is the label of an internal node and $b_{i+1}=\smac(b_i)$.
		\end{enumerate}
		Define the \emph{end of the maximal chain} of $t$ to be: $\emac(t):=b_k$.
		\item Let the unique path from the root to the leaf labeled $M$ be $m=c_1\to c_2\to\cdots\to c_l=M$. We call it the \emph{max-spine} of $t$ and define the \emph{parent of the maximum leaf} in $t$ to be: $\pom(t):=c_{l-1}$.
	\end{enumerate}
	\end{definition}
	By convention $\pom(t)=0$ and $\eoc(t)=\emac(t)=1$ for the unique $t\in\fT_1$. Even though the trees in $\fT_{2n+1}$ are nonplanar, it is to our convenience that we draw them as planar trees according to the following rule. For each internal node $a$ that is the parent of two subtrees $t_b$ and $t_c$ with $\max(t_b)<\max(t_c)$, we always draw $t_b$ as the left subtree and $t_c$ as the right subtree of $a$. A tree drawn this way is said to be \emph{standard}. The following observation will be useful.
	\begin{fact}\label{fact:left-right-branches}
	For a tree $t\in\fT_{2n+1}$ drawn in a standard way, the path from the root $1$ to $\emac(t)$ is exactly the \emph{left branch} of $t$ (i.e., the path that starts from the root and ends at a leaf, containing only left pointing edges), while the max-spine is exactly the \emph{right branch} of $t$ (i.e., the path that starts from the root and ends at a leaf, containing only right pointing edges), with $\pom(t)$ being the penultimate node on this branch.
	\end{fact}

	Table~\ref{tab:complete-5} lists all four trees in $\fT_5$ in a standard way, and their associated statistics are displayed in the middle three columns.

\begin{table}[htbp]
  \centering
  \caption{$\fT_5$ and $\AA_5$ with associated tree statistics}
	\label{tab:complete-5}
  \begin{tabular}{ccccc}
    \toprule
    $\fT_5$ & $\eoc$ & $\emac$ & $\pom$ & $\AA_5$ \\
    \midrule
    \bintree{
      \node (n1) at (0,0) {$1$};
      \node (n2) at (-1,1) {$2$};
      \node (n3) at (1,1) {$3$};
      \node (n4) at (0,2) {$4$};
      \node (n5) at (2,2) {$5$};
      \draw (n1) -- (n2) (n1) -- (n3) (n3) -- (n4) (n3) -- (n5);
    } & $2$ & $2$ & $3$ & $21435$ \\
    \addlinespace[6pt]
    \bintree{
      \node (n1) at (0,0) {$1$};
      \node (n3) at (-1,1) {$3$};
      \node (n2) at (1,1) {$2$};
      \node (n4) at (0,2) {$4$};
      \node (n5) at (2,2) {$5$};
      \draw (n1) -- (n3) (n1) -- (n2) (n2) -- (n4) (n2) -- (n5);
    } & $4$ & $3$ & $2$ & $31425$ \\
    \addlinespace[6pt]
    \bintree{
      \node (n1) at (0,0) {$1$};
      \node (n4) at (-1,1) {$4$};
      \node (n2) at (1,1) {$2$};
      \node (n3) at (0,2) {$3$};
      \node (n5) at (2,2) {$5$};
      \draw (n1) -- (n4) (n1) -- (n2) (n2) -- (n3) (n2) -- (n5);
    } & $3$ & $4$ & $2$ & $41325$ \\
    \addlinespace[6pt]
    \bintree{
      \node (n1) at (0,0) {$1$};
      \node (n2) at (-1,1) {$2$};
      \node (n5) at (1,1) {$5$};
      \node (n3) at (-2,2) {$3$};
      \node (n4) at (0,2) {$4$};
      \draw (n1) -- (n2) (n1) -- (n5) (n2) -- (n3) (n2) -- (n4);
    } & $3$ & $3$ & $1$ & $32415$ \\
    \bottomrule
  \end{tabular}
\end{table}

\emph{Andr\'e permutation (of the first kind)} was introduced by Foata and Sch\"utzenberger~\cite{FoataSchutzenberger1973} and it can be defined in several ways. For our purpose the following two versions are required.

\begin{definition}[recursive version]\label{def:And-1}
The empty word and each one-letter word are Andr\'e permutations. If $w=x_1x_2\cdots x_n$ ($n\ge 2$) is a permutation of a set of positive integers $Y=\{y_1<y_2<\cdots<y_n\}$, write $w=u\min(w)u'$. Then, $w$ is said to be an Andr\'e permutation if both $u$ and $u'$ are themselves Andr\'e permutations, and furthermore if $\max(uu')$ belongs to $u'$.
\end{definition}
It is not hard to deduce from the definition above that every Andr\'e permutation ends with its maximal letter. For a permutation $w=w_1\cdots w_n$, we write $w[j,k]:=w_jw_{j+1}\cdots w_k$ for the contiguous factor of $w$ from index $j$ to $k$.

\begin{definition}[interval version]\label{def:And-2}
A permutation $w=w_1w_2\cdots w_n$ with $w_n=n$ is said to be an Andr\'e permutation if there is no pair of indices $1\le j<k\le n-1$ with $k-j\ge 2$ such that
\begin{align}\label{Andre-2 condition}
\max(w[j,k])=w_j \text{ and } \min(w[j,k])=w_k.
\end{align}
\end{definition}

To make the paper self-contained, we include a proof of the equivalence between the two versions of definitions.
\begin{lemma}\label{lem:equiv}
The above two definitions of Andr\'e permutations are equivalent.
\end{lemma}
\begin{proof}
In the recursive decomposition of $w = w_1 \dots w_n$, we repeatedly split each contiguous factor $u$ at its minimum letter $w_i = \min(u)$ into $u = u' w_i u''$. Because earlier splits occur only at elements strictly smaller than $w_i$, the unique factor in the recursion whose minimum is $w_i$ corresponds to the \emph{maximal contiguous interval} $[\ell_i,r_i]$ containing $i$ such that $\min(w[\ell_i,r_i]) = w_i$\footnote{This is usually called the $x$-factorization of $w$ at $w_i$; see for instance \cite{FoataStrehl1974} and \cite[p.~9]{BrandenClaesson2011}.}. Thus $w$ is an Andr\'e permutation (recursive version) if and only if at every index $1 \le i \le n$, the factors $u' = w[\ell_i,i-1]$ and $u'' = w[i+1,r_i]$ satisfy either $u' = \emptyset$ or $\max(u') < \max(u'')$ (with $\max(\emptyset) := -\infty$).

\noindent\textbf{(Definition~\ref{def:And-1} $\implies$ Definition~\ref{def:And-2}).} Suppose $w=w_1\cdots w_n\in\fS_n$ and there exists a pair $(j, k)$ with $1 \le j < k \le n-1$ and $k - j \ge 2$ that satisfies~\eqref{Andre-2 condition}. Let $i \in \{j, \dots, k-1\}$ be the index where $w[j,k-1]$ achieves its minimum:
\[
w_i = \min(w[j,k-1]).
\]
Since $k - j \ge 2$ and $w_j = \max(w[j,k-1])$, we have $j < i \le k-1$. Consider the maximal interval $[\ell_i,r_i]$ containing $i$ with $\min(w[\ell_i, r_i]) = w_i$:
\begin{enumerate}
    \item Since $\min(w[j,i]) = w_i$, we have $\ell_i \le j < i$. Thus $w_j \in u' := w[\ell_i,i-1]$ and $\max(u') \ge w_j$.
    \item Since $w_k = \min(w[j,k]) < w_i$ and $k > i$, the interval $[\ell_i,r_i]$ cannot extend to $k$, so $r_i \le k-1$.
    \item Consequently, $[i+1,r_i] \subseteq [j+1,k-1]$, so every element of $u'':= w[i+1,r_i]$ is strictly smaller than $w_j = \max(w[j,k])$. Thus $\max(u'') < w_j \le \max(u')$, showing that $w$ is not Andr\'e in the sense of Definition~\ref{def:And-1}.
\end{enumerate}

\noindent\textbf{(Definition~\ref{def:And-2} $\implies$ Definition~\ref{def:And-1}).} Suppose $w$ is not Andr\'e (recursive version), then for some $i \in \{1, \dots, n\}$, the maximal interval $w[\ell_i,r_i]=u'w_iu''$ with minimum $w_i$ has $u' = w[\ell_i,i-1] \ne \emptyset$ and $\max(u') > \max(u'')$. Let $j \in [\ell_i, i-1]$ be the index where $u'$ achieves its maximum, thus $w_j = \max(w[\ell_i,r_i])$. 

Note that $r_i < n$: since $u' \ne \emptyset$, we have $w_i < w_j \le n$, so $i < n$; if $r_i = n$, then $w_n = n \in u''$, which would contradict $\max(u') > \max(u'')$. Because $w[\ell_i,r_i]$ is maximal on the right with respect to having elements no less than $w_i$, the element immediately to its right at index $k := r_i + 1$ must satisfy $w_k < w_i$. Moreover, $k \le n-1$ since $w_n = n > w_i$. We aim to show that the pair $(j,k)$ satisfies the condition \eqref{Andre-2 condition}, thereby implying that $w$ is not an Andr\'e permutation (interval version).
\begin{enumerate}
    \item Since $j < i < r_i + 1 = k$, we have $1 \le j < k \le n-1$ and $k - j \ge 2$.
    \item Since $[j, k-1] \subseteq [\ell_i, r_i]$, every element in $w[j,k-1]$ lies between $w_i$ and $w_j$.
    \item Because $w_k < w_i < w_j$, we obtain $\max(w[j,k]) = w_j$ and $\min(w[j,k]) = w_k$, so $(j, k)$ indeed satisfies~\eqref{Andre-2 condition}.
\end{enumerate}
\end{proof}

Recall that the \emph{reverse complement} is an involution $\rc:\fS_n\to \fS_n$ that sends every permutation $\pi=\pi_1\cdots\pi_n$ to $(n+1-\pi_n)(n+1-\pi_{n-1})\cdots(n+1-\pi_1)$.

\begin{lemma}\label{lem:trc}
The reverse complement induces an involution on the set of alternating Andr\'e permutations:
\begin{align*}
\trc: \AA_{2n+1} &\to \AA_{2n+1} \\
\pi_1\cdots\pi_{2n+1} &\mapsto (2n+1-\pi_{2n})(2n+1-\pi_{2n-1})\cdots(2n+1-\pi_1)\pi_{2n+1}.
\end{align*}
\end{lemma}
\begin{proof}
By Definition~\ref{def:And-2} (the interval version), we see that the factor $\pi[j,k]$ of $\pi$ satisfies condition \eqref{Andre-2 condition} in Definition~\ref{def:And-2} if and only if the factor $\trc(\pi)[2n-k+1,2n-j+1]$ of $\trc(\pi)$ satisfies this condition. Therefore, $\pi$ is an Andr\'e permutation if and only if $\trc(\pi)$ is an Andr\'e permutation. In addition, the mapping $\trc$ clearly sends an alternating permutation to an alternating permutation. Hence $\pi\in\AA_{2n+1}$ if and only if $\trc(\pi)\in\AA_{2n+1}$, meaning that $\trc$ is well-defined. It being an involution is clear from the definition.
\end{proof}

\begin{proof}[Proof of Theorem~\ref{thm:tree to perm}]
Given a tree $t\in\fT_{2n+1}$ drawn as standard, we read its labels by conducting an in-order traversal (i.e., left-root-right) of $t$, and let the resulting permutation be the image $\Psi(t)$. The relations in \eqref{stats:tree to perm} follow immediately from Fact~\ref{fact:left-right-branches}. The last column of Table~\ref{tab:complete-5} shows the permutations that correspond to the trees in the first column under the mapping $\Psi$. 

To see that $\Psi$ is well-defined, i.e., $\Psi(t)$ indeed belongs to $\AA_{2n+1}$, it suffices to note that $t$ being a strict tree ensures that $\Psi(t)$ is an alternating permutation, while drawing $t$ in a standard way forces $\Psi(t)$ to be an Andr\'e permutation (the recursive version). Conversely, for any $\pi\in\AA_{2n+1}$, we build the increasing binary tree of $\pi$ by following the decomposition $\pi = u \min(\pi) u'$ recursively. In this tree a letter $\pi_i$ has a left (resp.~right) child if and only if $\pi_{i-1} > \pi_i$ (resp.~$\pi_{i+1} > \pi_i$), with $\pi_0=\pi_{2n+2}=0$; since $\pi$ is down-up of odd length and ends with an ascent (recall that $\pi_{2n+1}=2n+1$), every letter is either a peak or a valley, so the tree is strict. The recursive Andr\'e condition (Definition~\ref{def:And-1}) says that at every node the right subtree contains the larger maximum, so the tree is standard. Let this tree be $\Psi^{-1}(\pi)$ and we see that indeed $\Psi^{-1}(\Psi(t))=t$ and $\Psi(\Psi^{-1}(\pi))=\pi$, making $\Psi$ a bijection.

Finally, the mapping $\trc$ constructed in Lemma~\ref{lem:trc} clearly implies \eqref{eq:reverse-complement}.
\end{proof}

\begin{proof}[Bijective proof of Theorem~\ref{thm:joint symm}]
Let $\Theta$ be the composed mapping
\[
	\Theta:=\Phi^{-1}\comp\Psi^{-1}\comp\trc\comp\Psi\comp\Phi: \fT_{2n+1}\to\fT_{2n+1}.
\]
Since both $\Phi$ and $\Psi$ are bijections while $\trc$ is an involution, it follows that $\Theta$ is also an involution. Making use of the relations \eqref{pom=pom}, \eqref{emac=eoc}, \eqref{stats:tree to perm}, and the definition of $\trc$, we can compute that for each $t\in\fT_{2n+1}$,
\begin{align*}
\eoc(\Theta(t)) &= \eoc(\Phi^{-1}(\Psi^{-1}(\trc(\Psi(\Phi(t))))))\\
&= \emac(\Psi^{-1}(\trc(\Psi(\Phi(t)))))\\
&= \trc(\Psi(\Phi(t)))_1\\
&= 2n+1-\Psi(\Phi(t))_{2n}\\
&= 2n+1-\pom(\Phi(t))\\
&= 2n+1-\pom(t).
\end{align*}
A similar computation verifies that $\pom(\Theta(t))=2n+1-\eoc(t)$ and completes the proof.
\end{proof}

\section{Proof of Theorem~\ref{thm:eoc to emac}}\label{sec:eoc-emac}
\begin{figure}
\begin{center}
\resizebox{0.9\textwidth}{!}{%
\begin{tikzpicture}[
    scale=1.0, 
    every node/.style={font=\small, inner sep=2pt}
]
    \begin{scope}[shift={(0,0)}]
        \node (L_p1)    at (0, 0)       {$p_1$};
        \node (L_p2)    at (1.2, 1.0)   {$p_2$};
        \node (L_dots1) at (2.2, 1.85)  {\rotatebox{90}{$\ddots$}};
        \node (L_pis)   at (3.2, 2.7)   {$p_{i^*}$};
        \node (L_dots2) at (4.2, 3.55)  {\rotatebox{90}{$\ddots$}};
        \node (L_ps)    at (5.2, 4.4)   {$p_s$};
        \node (L_ps+1)  at (6.2, 5.3)   {$p_{s+1}$};

        \node (L_s1)    at (-0.7, 1.1)  {$S_1$};
        \node (L_s2)    at (0.5, 2.1)   {$S_2$};
        \node (L_sis)   at (2.5, 3.8)   {$S_{i^*}$};
        \node (L_ss)    at (4.5, 5.5)   {$S_s$};

        \draw (L_p1) -- (L_p2);
        \draw (L_p2) -- (L_dots1);
        \draw (L_dots1) -- (L_pis);
        \draw (L_pis) -- (L_dots2);
        \draw (L_dots2) -- (L_ps);

        \draw (L_p1) -- (L_s1);
        \draw (L_p2) -- (L_s2);
        \draw (L_pis) -- (L_sis);
        \draw (L_ps) -- (L_ss);
        \draw (L_ps) -- (L_ps+1);
    \end{scope}

    \draw[->, thick] (7.3, 3) -- node[above, font=\normalsize] {$\Phi$} (8.7, 3);

    \begin{scope}[shift={(10.0,0)}]
        \node (R_p1)    at (0, 0)       {$p_1$};
        \node (R_p2)    at (1.2, 1.0)   {$p_2$};
        \node (R_dots1) at (2.2, 1.85)  {\rotatebox{90}{$\ddots$}};
        \node (R_pis)   at (3.2, 2.7)   {$p_{i^*}$};
        \node (R_dots2) at (4.2, 3.55)  {\rotatebox{90}{$\ddots$}};
        \node (R_ps)    at (5.2, 4.4)   {$p_s$};
        \node (R_ps+1)  at (6.2, 5.3)   {$p_{s+1}$};

        \node (R_s1)    at (-0.85, 1.15) {$\Phi(S_{i^*})$};
        \node (R_s2)    at (0.5, 2.1)    {$S_1$};
        \node (R_sis)   at (2.5, 3.8)    {$S_{i^*-1}$};
        \node (R_ss)    at (4.5, 5.5)    {$S_s$};

        \draw (R_p1) -- (R_p2);
        \draw (R_p2) -- (R_dots1);
        \draw (R_dots1) -- (R_pis);
        \draw (R_pis) -- (R_dots2);
        \draw (R_dots2) -- (R_ps);

        \draw (R_p1) -- (R_s1);
        \draw (R_p2) -- (R_s2);
        \draw (R_pis) -- (R_sis);
        \draw (R_ps) -- (R_ss);
        \draw (R_ps) -- (R_ps+1);
    \end{scope}
\end{tikzpicture}%
}
\end{center}
\caption{The bijection $\Phi$.}\label{fig:Phi}
\end{figure}
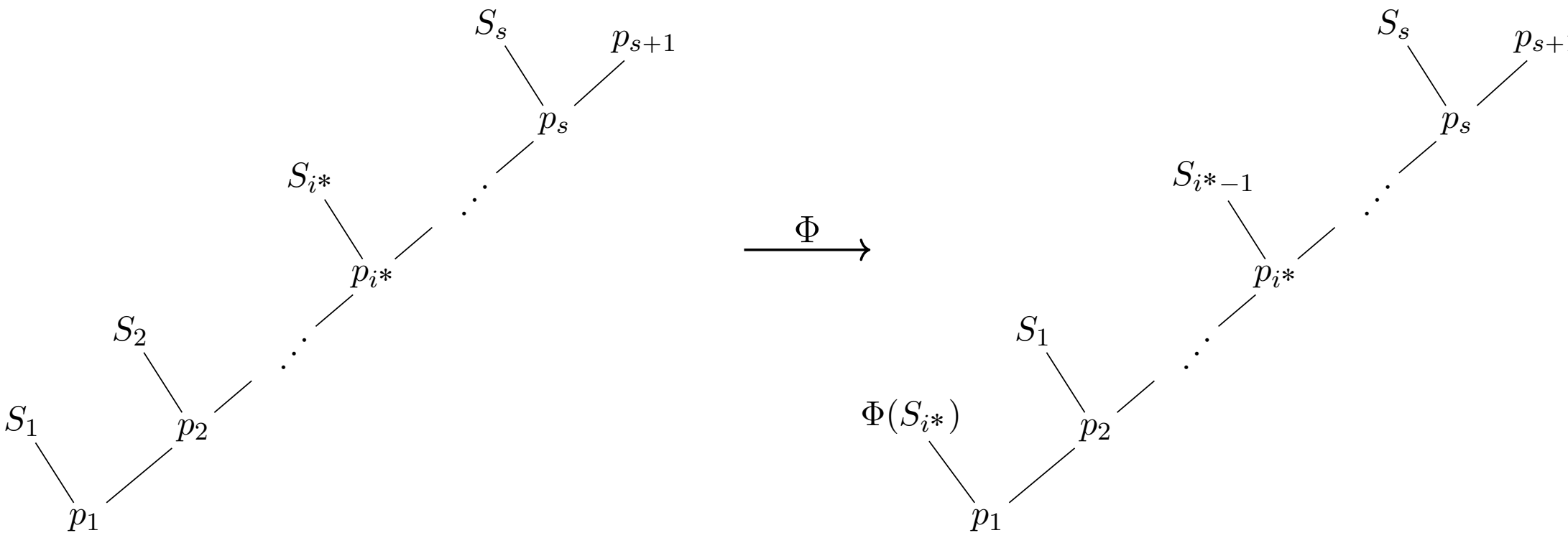

It remains to prove Theorem~\ref{thm:eoc to emac}. The first step is to understand how $\eoc$ and $\emac$ branch from the max-spine. For $n\ge 1$ and a tree $t\in\fT_{2n+1}$, suppose its max-spine is
\[
1 = p_1 \to p_2 \to \dots \to p_s \to p_{s+1} = 2n+1 \qquad (s \ge 1).
\]
At each internal node $p_i$ ($1 \le i \le s$) along the $\max$-spine, one of its two children is the next spine vertex $p_{i+1}$, and its other child is the root $r_i := \min(S_i)$ of an \emph{off-spine subtree} $S_i$. We write
\[
	t=(p_1,\ldots,p_s,p_{s+1};S_1,\ldots,S_s).
\]
An illustration is given in the left tree of Fig.~\ref{fig:Phi}.

We already observed in Fact~\ref{fact:left-right-branches} that $\emac(t)$ is on the end of the left branch of $t$, so in particular
\begin{align}\label{eq:emac_branch}
	\emac(t)=\emac(S_1).
\end{align}
On the other hand, let $i^* \in \{1, 2, \dots, s\}$ be the \emph{smallest} index on the max-spine such that
\begin{equation}\label{eq:istar_def}
r_{i^*} = \min(S_{i^*}) < p_{i^*+1}.
\end{equation}
(Such an index $i^*$ always exists because at $i = s$, $p_{s+1} = 2n+1 > r_s$.) Then the minimal chain of $t$ follows the max-spine from $p_1$ up to $p_{i^*}$, branches off at $p_{i^*}$ into $S_{i^*}$, and terminates at
\begin{equation}\label{eq:eoc_branch}
\eoc(t) = \eoc(S_{i^*}).
\end{equation}

\begin{proof}[Proof of Theorem~\ref{thm:eoc to emac}]
We are going to construct the bijection $\Phi$ recursively. Since both $\fT_1$ and $\fT_3$ contain a single tree, $\Phi$ acts as identity on them. Now for a given $n\ge 2$, assume by induction that the bijection $\Phi:\fT_{2m+1}\to\fT_{2m+1}$ satisfying \eqref{pom=pom} and \eqref{emac=eoc} has already been defined\footnote{Note that we tacitly extend the domains of both $\Phi$ and $\Phi^{-1}$ to strict trees on any finite label set $L\subset\N$ (either directly, since the construction never uses the actual labels, or by order-preserving relabeling), and that $\Phi$ preserves the label set, hence the root.} for $1\le m<n$. Take a tree $t=(p_1,\ldots,p_s,p_{s+1};S_1,\ldots,S_s)\in\fT_{2n+1}$, its image $\Phi(t)$ is constructed as follows; an illustration of $\Phi(t)$ can be found in Fig.~\ref{fig:Phi} (the tree on the right).
\begin{itemize}
	\item Find the smallest index $i^*$ as defined in \eqref{eq:istar_def}. 
	\item Cyclically shift the off-spine subtrees $S_1,\ldots,S_{i^*}$ upwards away from the root $p_1$. I.e., for $1\le j\le i^*-1$, detach the subtree $S_j$ from $p_j$ and attach it to $p_{j+1}$ as a subtree; detach $S_{i^*}$ from $p_{i^*}$ and attach $\Phi(S_{i^*})$ to $p_{1}$ as its new off-spine subtree. Let this newly obtained tree be 
	\[
		\Phi(t) := (p_1,\ldots,p_s,p_{s+1};\Phi(S_{i^*}),S_1,\ldots,S_{i^*-1},S_{i^*+1},\ldots,S_s).
	\]
\end{itemize}
To finish the proof, we need to show the following items.
\begin{enumerate}
	\item $\Phi(t)$ indeed belongs to $\fT_{2n+1}$, thereby showing that $\Phi$ is well-defined.
	\item The statistics are transformed according to \eqref{pom=pom} and \eqref{emac=eoc}.
	\item $\Phi$ is invertible.
\end{enumerate}

To see (1), it is clear that $\Phi(t)$ as drawn in Fig.~\ref{fig:Phi} continues to be strict. It remains to verify that $\Phi(t)$ is still increasingly labeled. Since $p_1=\min(\Phi(t))$, attaching $\Phi(S_{i^*})$ to it is allowed; for every $j\in\{1,\ldots,i^*-1\}$, the minimality of the index $i^*$ (see \eqref{eq:istar_def}) ensures that $r_j=\min(S_j)>p_{j+1}$, so attaching $S_j$ to $p_{j+1}$ is also allowed.

For (2), since preserving max-spine is actually stronger than preserving $\pom$ (the penultimate node $p_s$ in the max-spine), we have \eqref{pom=pom}. Using both \eqref{eq:emac_branch} and \eqref{eq:eoc_branch}, we have
\[
	\emac(\Phi(t))=\emac(\Phi(S_{i^*}))=\eoc(S_{i^*})=\eoc(t),
\]
where we need to use the induction hypothesis for the middle ``$=$''. This gives us \eqref{emac=eoc}.

Finally for (3), suppose we are given a tree $t=(p_1,\ldots,p_s,p_{s+1};R_1,\ldots,R_s)\in\fT_{2n+1}$ with $r_i:=\min(R_i)$, $1\le i\le s$, we take two steps to recover the preimage tree $\Phi^{-1}(t)$. 
\begin{itemize}
	\item Find the unique index $j^*\in\{1,\ldots,s\}$, such that
	\begin{align}\label{def:jstar}
		p_{j^*}<r_1<p_{j^*+1}.
	\end{align}
	\item Detach $R_1$ from $p_1$ and attach $\Phi^{-1}(R_1)$ to $p_{j^*}$; cyclically shift the off-spine subtrees $R_2,\ldots,R_{j^*}$ down towards the root $p_1$, i.e., for $2\le i\le j^*$, detach the subtree $R_i$ from $p_i$ and attach it to $p_{i-1}$ as a subtree. Let this newly obtained tree be
	\[
		\Phi^{-1}(t) := (p_1,\ldots,p_s,p_{s+1};R_2,\ldots,R_{j^*},\Phi^{-1}(R_1),R_{j^*+1},\ldots,R_s).
	\]
\end{itemize}
Note that for $2\le i\le j^*$, attaching $R_i$ to $p_{i-1}$ is always allowed since
\begin{align}\label{shift-down valid}
	r_i=\min(R_i)>p_i>p_{i-1};
\end{align}
condition~\eqref{def:jstar} ensures the validity of attaching $\Phi^{-1}(R_1)$ to $p_{j^*}$. Moreover, combining \eqref{def:jstar} and \eqref{shift-down valid} we deduce that $j^*$ coincides with the smallest index $i^*$ for $\Phi^{-1}(t)$ as defined by \eqref{eq:istar_def}. Consequently, we see that $\Phi^{-1}(t)$ indeed belongs to $\fT_{2n+1}$, and that
\[
	\Phi(\Phi^{-1}(t))=t=\Phi^{-1}(\Phi(t)).
\]
So $\Phi$ is invertible and the proof of Theorem~\ref{thm:eoc to emac} is now complete.
\end{proof}


To illustrate the recursive nature of the transformation $\Phi$ and the composed involution $\Theta$ when multiple nested cyclic shifts occur, let $n=7$ and start with the tree $t\in\fT_{15}$ given below. We show step-by-step how to compute its image $\Theta(t)$, with illustrations of all intermediate trees included at the end. Throughout, a tree consisting of a single vertex $a$ is simply denoted by $(a)$. Let
\begin{align*}
	&t:=(1,2,3,11,15;S_1,S_2,S_3,S_4), \text{ where } S_1=(13),~S_2=(14),\\
	&S_3=(4,5,10;U_1,U_2),~S_4=(12), \text{ and }U_1=(9),~U_2=(6,8;(7)).
\end{align*}
Observe that $\eoc(t)=7$ and $\pom(t)=11$.

\begin{figure}[htbp]
\centering
\begin{tikzpicture}[
    scale=0.72,
    every node/.style={font=\small},
    v/.style={inner sep=2pt},
    edge/.style={thick},
    maparrow/.style={-{Stealth[length=6pt, width=5pt]}, thick}
]

\begin{scope}[shift={(0,0)}]
    \node[v] (TL1)  at (0, 0)     {$1$};
    \node[v] (TL13) at (-0.8, 0.9){$13$};
    \node[v] (TL2)  at (0.9, 0.9) {$2$};
    \node[v] (TL14) at (0.1, 1.8) {$14$};
    \node[v] (TL3)  at (1.8, 1.8) {$3$};
    \node[v] (TL4)  at (0.8, 2.7) {$4$};
    \node[v] (TL11) at (2.8, 2.7) {$11$};
    \node[v] (TL12) at (2.1, 3.6) {$12$};
    \node[v] (TL15) at (3.6, 3.6) {$15$};

    \node[v] (TL9)  at (-0.2, 3.6){$9$};
    \node[v] (TL5)  at (1.3, 3.8) {$5$};
    \node[v] (TL6)  at (0.5, 4.7) {$6$};
    \node[v] (TL7)  at (-0.3, 5.6){$7$};
    \node[v] (TL8)  at (1.2, 5.6) {$8$};
    \node[v] (TL10) at (2.1, 4.9) {$10$};

    \draw[edge] (TL1)--(TL13) (TL1)--(TL2) (TL2)--(TL14) (TL2)--(TL3)
                (TL3)--(TL4)  (TL3)--(TL11) (TL11)--(TL12) (TL11)--(TL15)
                (TL4)--(TL9)  (TL4)--(TL5) (TL5)--(TL6) (TL5)--(TL10)
                (TL6)--(TL7)  (TL6)--(TL8);

    \draw[dashed, thick, rounded corners=14pt]
        (0.8, 2.2) -- (-0.85, 3.6) -- (-0.85, 6.1) -- (2.85, 5.65) -- (1.65, 3.55) -- cycle;
    \node at (-1.5, 4.5) {$S_3=$};

    \node[below=0.45cm of TL1, xshift=1.8cm, font=\normalsize] (TL_cap)
        {$t \in \fT_{15}\ (\mathrm{eoc} = 7,\ \mathrm{pom} = 11)$};
\end{scope}

\begin{scope}[shift={(10.8,0)}]
    \node[v] (TR1)  at (0, 0)     {$1$};
    \node[v] (TR4)  at (-1.4, 1.0){$4$};
    \node[v] (TR6)  at (-2.3, 2.1){$6$};
    \node[v] (TR7)  at (-3.0, 3.1){$7$};
    \node[v] (TR8)  at (-1.7, 3.1){$8$};
    \node[v] (TR5)  at (-0.6, 2.1){$5$};
    \node[v] (TR9)  at (-1.1, 3.1){$9$};
    \node[v] (TR10) at (-0.1, 3.1){$10$};

    \node[v] (TR2)  at (1.4, 1.0) {$2$};
    \node[v] (TR13) at (0.7, 2.0) {$13$};
    \node[v] (TR3)  at (2.3, 2.0) {$3$};
    \node[v] (TR14) at (1.6, 3.0) {$14$};
    \node[v] (TR11) at (3.2, 3.0) {$11$};
    \node[v] (TR12) at (2.5, 4.0) {$12$};
    \node[v] (TR15) at (4.0, 4.0) {$15$};

    \draw[edge] (TR1)--(TR4) (TR4)--(TR6) (TR6)--(TR7) (TR6)--(TR8)
                (TR4)--(TR5) (TR5)--(TR9) (TR5)--(TR10)
                (TR1)--(TR2) (TR2)--(TR13) (TR2)--(TR3) (TR3)--(TR14)
                (TR3)--(TR11) (TR11)--(TR12) (TR11)--(TR15);

    \draw[dashed, thick, rounded corners=14pt]
        (-1.4, 0.5) -- (-3.8, 3.6) -- (.8, 3.6) -- cycle;
    \node at (-4, 2) {$\Phi(S_3)=$};

    \node[below=0.45cm of TR1, xshift=0.5cm, font=\normalsize] (TR_cap)
        {$\Phi(t) \in \fT_{15}\ (\emac = 7,\ \mathrm{pom} = 11)$};
\end{scope}

\begin{scope}[shift={(0,-9.2)}]
    \node[v] (BL1)  at (0, 0)     {$1$};
    \node[v] (BL13) at (-0.8, 0.9){$13$};
    \node[v] (BL2)  at (0.9, 0.9) {$2$};
    \node[v] (BL3)  at (0.0, 1.9) {$3$};
    \node[v] (BL4)  at (-0.6, 2.9){$4$};
    \node[v] (BL12) at (0.6, 2.9) {$12$};
    \node[v] (BL5)  at (1.9, 1.8) {$5$};
    \node[v] (BL14) at (1.3, 2.8) {$14$};
    \node[v] (BL6)  at (2.8, 2.7) {$6$};
    \node[v] (BL10) at (2.2, 3.7) {$10$};
    \node[v] (BL7)  at (3.7, 3.6) {$7$};
    \node[v] (BL11) at (3.1, 4.6) {$11$};
    \node[v] (BL8)  at (4.6, 4.5) {$8$};
    \node[v] (BL9)  at (4.0, 5.5) {$9$};
    \node[v] (BL15) at (5.5, 5.4) {$15$};

    \draw[edge] (BL1)--(BL13) (BL1)--(BL2) (BL2)--(BL3) (BL3)--(BL4) (BL3)--(BL12)
                (BL2)--(BL5)  (BL5)--(BL14) (BL5)--(BL6) (BL6)--(BL10)
                (BL6)--(BL7)  (BL7)--(BL11) (BL7)--(BL8) (BL8)--(BL9) (BL8)--(BL15);

    \node[below=0.45cm of BL1, xshift=2.2cm, font=\normalsize] (BL_cap)
        {$\Theta(t) \in \fT_{15}\ (\mathrm{eoc} = 4,\ \mathrm{pom} = 8)$};
\end{scope}

\begin{scope}[shift={(10.8,-9.2)}]
    \node[v] (BR1)  at (0, 0)     {$1$};
    \node[v] (BR3)  at (-1.1, 1.0){$3$};
    \node[v] (BR4)  at (-1.8, 2.0){$4$};
    \node[v] (BR12) at (-0.4, 2.0){$12$};

    \node[v] (BR2)  at (1.1, 1.0) {$2$};
    \node[v] (BR13) at (0.5, 2.0) {$13$};
    \node[v] (BR5)  at (2.0, 1.9) {$5$};
    \node[v] (BR14) at (1.4, 2.9) {$14$};
    \node[v] (BR6)  at (2.9, 2.8) {$6$};
    \node[v] (BR10) at (2.3, 3.8) {$10$};
    \node[v] (BR7)  at (3.8, 3.7) {$7$};
    \node[v] (BR11) at (3.2, 4.7) {$11$};
    \node[v] (BR8)  at (4.7, 4.6) {$8$};
    \node[v] (BR9)  at (4.1, 5.6) {$9$};
    \node[v] (BR15) at (5.6, 5.5) {$15$};

    \draw[edge] (BR1)--(BR3) (BR3)--(BR4) (BR3)--(BR12)
                (BR1)--(BR2) (BR2)--(BR13) (BR2)--(BR5) (BR5)--(BR14)
                (BR5)--(BR6) (BR6)--(BR10) (BR6)--(BR7) (BR7)--(BR11)
                (BR7)--(BR8) (BR8)--(BR9) (BR8)--(BR15);

    \node[below=0.45cm of BR1, xshift=1.8cm, font=\normalsize] (BR_cap)
        {$t^* \in \fT_{15}\ (\emac = 4,\ \mathrm{pom} = 8)$};
\end{scope}

\draw[maparrow] (4.8, 2.6) -- node[above=4pt, font=\normalsize] {$\Phi$} (6.8, 2.6);

\draw[maparrow] (7.4, -6.4) -- node[above=4pt, font=\normalsize] {$\Phi^{-1}$} (5.4, -6.4);

\draw[maparrow] (11.5, -1.8) -- node[right=6pt, font=\normalsize] {$\Psi^{-1} \comp \widetilde{\mathrm{rc}} \comp \Psi$} (11.5, -3.5);

\draw[maparrow, dashed] (1.6, -1.8) -- node[right=6pt, font=\normalsize] {$\Theta = \Phi^{-1} \comp \Psi^{-1} \comp \widetilde{\mathrm{rc}} \comp \Psi \comp \Phi$} (1.6, -3.5);

\end{tikzpicture}
\caption{The composed involution $\Theta$ on $\fT_{15}$.}
\label{fig:tree-mapping-diagram}
\end{figure}
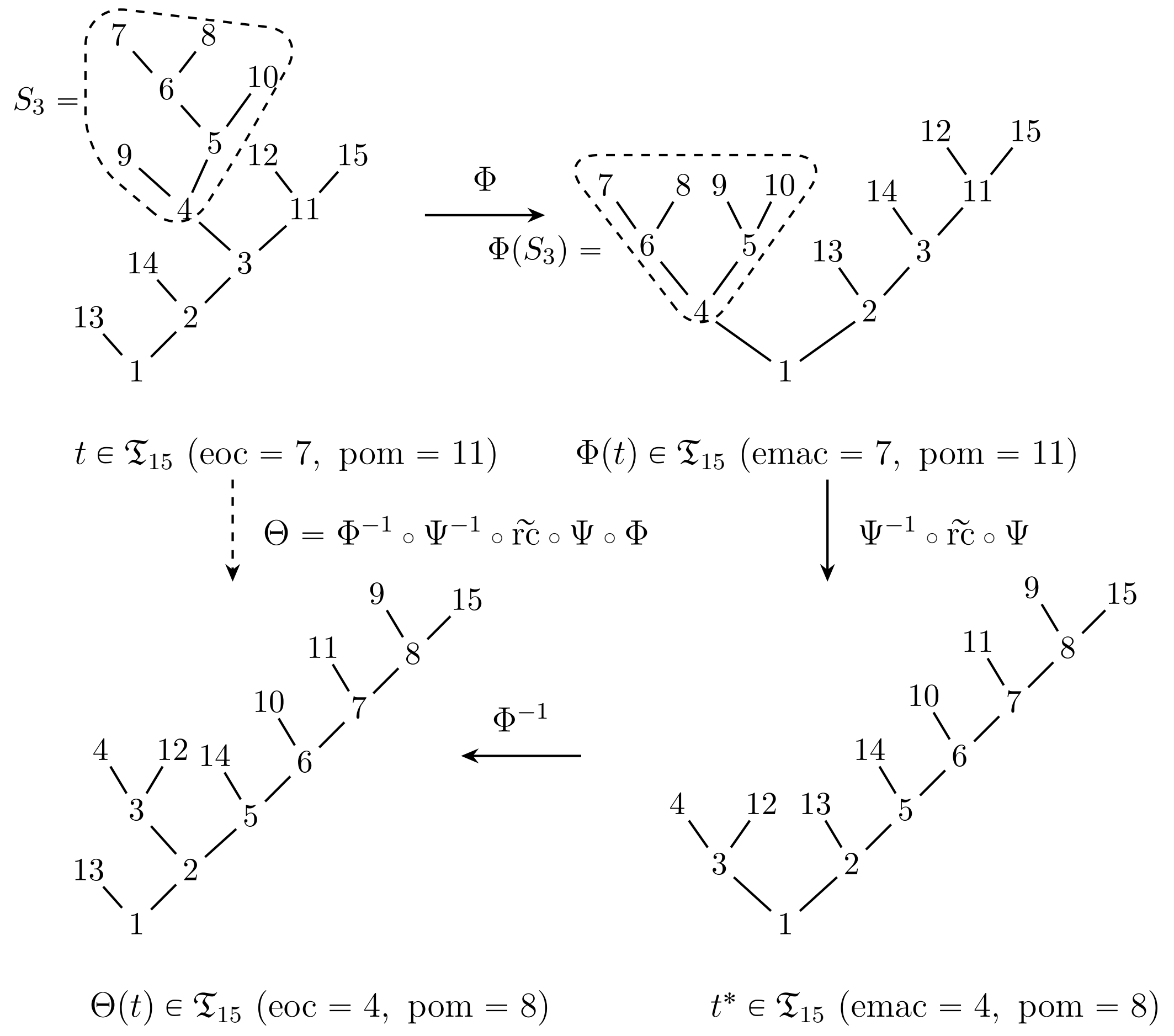

\noindent {\bf Step 1: find $\Phi(t)$ via three nested cyclic shiftings.}

The index $i^*$ for $t$ is seen to be $i^*=3$, so $\Phi(t)=(1,2,3,11,15;\Phi(S_3),S_1,S_2,S_4)$. To find $\Phi(S_3)$, we find the critical index $i^*$ for the subtree $S_3$ to be $i^*=2$, thus
\[
	\Phi(S_3)=(4,5,10;\Phi(U_2),U_1)=(4,5,10;U_2,U_1),
\]
since $\Phi$ acts on the three-vertices tree $U_2$ as identity.

\noindent {\bf Step 2: applying the composed mapping $\Psi^{-1}\comp\trc\comp\Psi$.}
An in-order traverse of the tree $\Phi(t)$ renders the alternating Andr\'e permutation
\[
	\pi:=\Psi(\Phi(t)) = \mathbf{7}\;\; 6\;\; 8\;\; 4\;\; 9\;\; 5\;\; 10\;\; 1\;\; 13\;\; 2\;\; 14\;\; 3\;\; 12\;\; \mathbf{11}\;\; 15\in \AA_{15}, 
\]
where $\pi_1=7=\eoc(t)$ and $\pi_{14}=11=\pom(t)$. Applying $\trc$ we get
\[
	\pi^*:=\trc(\pi)=\mathbf{4}\;\; 3\;\; 12\;\; 1\;\; 13\;\; 2\;\; 14\;\; 5\;\; 10\;\; 6\;\; 11\;\; 7\;\; 9\;\; \mathbf{8}\;\; 15 \in \AA_{15}.
\]
Observe that $\pi^*_1=4=15-\pom(t)$ and $\pi^*_{14}=8=15-\eoc(t)$. Applying $\Psi^{-1}$ we recover the tree
\[
	t^*:=\Psi^{-1}(\pi^*)=(1,2,5,6,7,8,15;R_1,R_2,R_3,R_4,R_5,R_6),
\]
where $R_1=(3,12;(4))$, $R_2=(13)$, $R_3=(14)$, $R_4=(10)$, $R_5=(11)$, and $R_6=(9)$.

\noindent {\bf Step 3: find the final image $\Theta(t)$ via $\Phi^{-1}$.}
To deduce $\Phi^{-1}(t^*)$, we check the first off-spine subtree $R_1=(3,12;(4))$, whose root is $3$, implying that the unique index $j^*=2$, i.e., we should cyclically shift $R_1$ and $R_2$. Note that $\Phi^{-1}(R_1)=R_1$. Consequently, we have found our final image tree to be
\[
	\Theta(t)=\Phi^{-1}(t^*)=(1,2,5,6,7,8,15;R_2,R_1,R_3,R_4,R_5,R_6).
\]
One checks that $\eoc(\Theta(t))=4=15-\pom(t)$ and $\pom(\Theta(t))=8=15-\eoc(t)$, as expected.

\section*{Acknowledgement} 
The author is grateful to Kathy Q. Ji for organizing a stimulating conference in February 2026. During this gathering, Guo-Niu Han kindly shared his interest in a simpler proof of Theorem~\ref{thm:joint symm}, explained the subtleties and difficulties of their original argument, and pointed out several relevant references, including \cite{Han2013}. The author further thanks him for his insightful comments on an earlier version of this manuscript.



\begin{thebibliography}{99}

		\bibitem{Andre1879}
		D. Andr\'e, D\'eveloppement de $\sec x$ and $\tan x$, \emph{C. R. Acad. Sci. Paris} {\bf88} (1879), 965--967.

		\bibitem{BrandenClaesson2011}
		P. Br\"and\'en and A. Claesson, \emph{Mesh patterns and the expansion of permutation statistics as sums of permutation patterns}, Electronic Journal of Combinatorics \textbf{18}(2) (2011), \#P5.

		\bibitem{Comtet1974}
		L. Comtet, \emph{Advanced Combinatorics}, D. Reidel, Dordrecht-Holland, Boston, MA, 1974.



		\bibitem{FoataHan2014}
		D. Foata and G.-N. Han, \emph{Tree calculus for bivariate difference equations}, Journal of Difference Equations and Applications \textbf{20}:11 (2014), 1453--1488.

		\bibitem{FoataSchutzenberger1973}
		D. Foata and M.-P. Sch\"utzenberger, \emph{Nombres d'Euler et permutations alternantes}, in \emph{A Survey of Combinatorial Theory}, North-Holland, Amsterdam, 1973, pp.~173--187.

		\bibitem{FoataStrehl1974}
		D. Foata and V. Strehl, \emph{Rearrangements of the symmetric group and enumerative properties of the tangent and secant numbers}, Mathematische Zeitschrift \textbf{137} (1974), 257--264.

		\bibitem{Han2013}
		G.-N. Han, \emph{The Poupard statistics on tangent and secant trees}, Strasbourg, preprint (2013).

		\bibitem{Poupard1989}
		C. Poupard, \emph{Deux propri\'et\'es des arbres binaires ordonn\'es stricts}, European Journal of Combinatorics \textbf{10} (1989), 369--374.

		\bibitem{OEIS} N.J.A. Sloane, On-Line Encyclopedia of Integer Sequences, available at \url{http://oeis.org}.

		\bibitem{Stanley2010}
		R.P. Stanley, \emph{A survey of alternating permutations}, Contemp. Math. {\bf531} (2010), 165--196.
	
		
	\end{thebibliography}
\end{document}